\documentclass[final,5p,twocolumn]{elsarticle}
\usepackage{lineno,hyperref}

\usepackage{amsmath}
\usepackage{comment}
\usepackage{setspace}
\usepackage{float}
\usepackage{color}
\usepackage{ulem}
\usepackage{mathptmx}

\modulolinenumbers[5]
\biboptions{numbers,sort&compress}

\graphicspath{{./figures/}}

\begin{document}

\begin{frontmatter}

\title{Novel flow field design for vanadium redox flow batteries via topology optimization}

\author[mymainaddress]{Chih-Hsiang Chen}
\author[mymainaddress]{Kentaro Yaji\corref{mycorrespondingauthor}}
\ead{yaji@mech.eng.osaka-u.ac.jp}
\author[mymainaddress]{Shintaro Yamasaki}
\author[mymainaddress]{Shohji Tsushima}
\author[mymainaddress]{Kikuo Fujita}


\cortext[mycorrespondingauthor]{Corresponding author}

\address[mymainaddress]{Department of Mechanical Engineering, Graduate School of Engineering, Osaka University, 2-1, Yamadaoka, Suita, Osaka 565-0871, Japan}

\begin{abstract}
This paper presents a three-dimensional topology optimization method for the design of flow field in vanadium redox flow batteries (VRFBs).
We focus on generating a novel flow field configuration for VRFBs via topology optimization, which has been attracted attention as a powerful design tool based on numerical optimization. An attractive feature of topology optimization is that a topology optimized configuration can be automatically generated without presetting a promising design candidate. 
In this paper, we formulate the topology optimization problem as a maximization problem of the electrode surface concentration in the negative electrode during the charging process.
The aim of this optimization problem is to obtain a topology optimized flow field that enables the improvement of mass transfer effect in a VRFB. 
We demonstrate that a novel flow field configuration can be obtained  through the numerical investigation.
To clarify the performance of the topology optimized flow field, we investigate the mass transfer effect through the comparison with reference flow fields---parallel and interdigitated flow fields---and the topology optimized flow field. 
In addition, we discuss the power loss that takes account of the polarization loss and pumping power, at various operating conditions.
\end{abstract}

\begin{keyword}
Redox flow battery\sep Flow field design \sep Topology optimization\sep Mass transfer effect
\end{keyword}

\end{frontmatter}


\section{Introduction}

\begin{figure*}[t]
	\centering
\includegraphics[width=150mm]{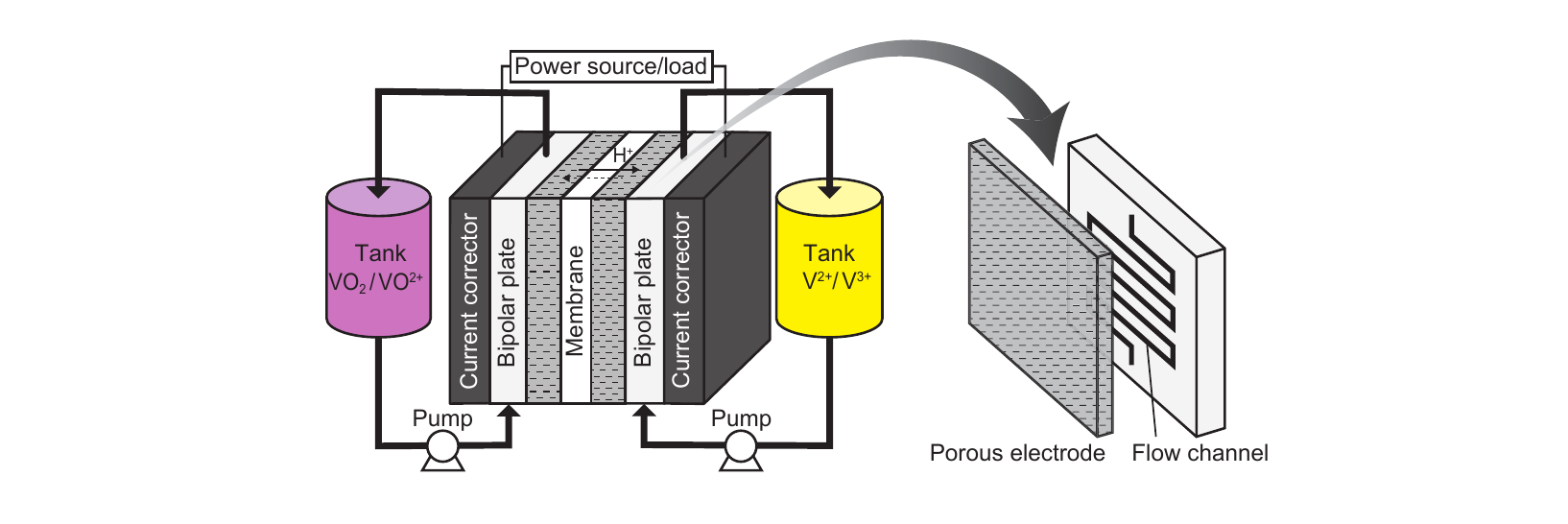}
	\vspace{-5mm}
	\caption{Schematic diagram of a redox flow battery.}
	\label{fig: sec1_1_1}
\end{figure*}

Over the past years, there has been growing interest in renewable energy since the increase of carbon emission poses a great threat to the environment. Despite the significance of renewable energy, the intermittent characteristics of sources, such as solar, wind or water, is a fatal drawback for renewable energy, which leads to the increased uncertainty in the supply of electricity. To address this issue, one promising solution is to regulate the power delivery via energy storage technology. Among the energy storage systems, vanadium redox flow batteries (VRFBs) attract a lot of attention due to the advantageous features:  scalability, low cost and long cycle life \cite{alotto2014redox}. However, achieving high performance in terms of power density is a critical issue for cost-effectiveness of VRFBs.

The polarization losses in VRFBs are mainly caused by ohmic, mass transfer and charge transfer losses \cite{milshtein2017quantifying}. Several researchers contributed to reduce overpotentials via a new cell architecture  \cite{aaron2012dramatic}, and modified electrode configuration and membrane \cite{sun1992modification,wang2007investigation,chen2013optimizing}. For charge transfer losses, Li et al.~\cite{li2013bismuth} proposed electrodes containing nanoparticles that can improve performance of VRFBs due to the acceleration of charge transfer compared to the conventional VRFBs. Li et al.~\cite{li2011graphite} claimed that reduction of charge transfer losses can be achieved by adding graphite oxide to electrodes. 

For given electrochemical conditions, the performance of VRFBs depends on mass transfer losses mainly, where mass transfer effect can be ameliorated by different flow fields \cite{zhou2017critical}. 
Xu et al.~\cite{xu2013numerical} numerically investigated the performance of VRFBs with several different types of flow fields and found the VRFB with serpentine flow field achieved maximum power-based efficiency at the optimal flow rate. Studies based on experiments demonstrated that interdigitated flow fields can further improve the performance of VRFBs due to the enhanced mass transfer effect \cite{tsushima2014efficient,darling2014influence,houser2016influence}. Although few types of flow field have been found to reduce mass transfer losses, it is still laborious for the design of flow field due to the complicated physical and chemical mechanisms in VRFBs.

Based on physical principles and mathematical models, Bends{\o}e and Kikuchi \cite{bendsoe1988generating} proposed topology optimization that is a powerful approach to find optimal configurations. 
Topology optimization expresses a structural optimization problem as a material distribution problem in a given design domain and then derives promising configuration on the basis of mathematical programming \cite{bendsoe2003topology}. 
One of the attractive feature of topology optimization is that an innovative configuration can be automatically generated from a blank design domain without designer's intuition.
Due to its high degree of design freedom, topology optimization has been applied to various structural optimization problems, e.g., stiffness maximization problems \cite{bendsoe1988generating,bendsoe1989optimal}, eigenfrequency problems \cite{diaz1992solutions,ma1995topological}, thermal problems \cite{li1999shape,iga2009topology}, and electromagnetic problems \cite{nomura2007structural,yamasaki2011level}; furthermore, applications to practical device designs have also been attracted attention in micro actuator design \cite{sigmund2001design}, fuel cell design \cite{iwai2011power,song20132d} and so on.

For flow field design problems, Borrvall and Petersson \cite{borrvall2003topology} proposed a topology optimization method to minimize power dissipation in Stokes flow, and this has been expanded to laminar Navier-Stokes flow problems \cite{gersborg2005topology,olesen2006high,kubo2017level} and turbulence problems \cite{yoon2016topology,dilgen2018topology}.
The fluid topology optimization has been applied to multiphysics problems such as fluid-structure interaction problems \cite{yoon2010topology,jenkins2015level}, forced convection problems \cite{matsumori2013topology,yaji2015topology,yaji2018large}, natural convection problems \cite{alexandersen2014topology,coffin2016level,alexandersen2016large} and  turbulent heat transfer problems \cite{kontoleontos2013adjoint,dilgen2018density}.

Recently, 
Yaji et al. \cite{yaji2018topology} proposed a topology optimization method for the design of flow fields in VRFBs. 
In their approach, instead of the formula of practical electrochemical reactions, a simplified formula is introduced as a two-dimensional model.
They provided novel flow field configurations of a VRFB and clarified that the optimized configurations tend to be the type of the interdigitated flow field.

As a more comprehensive study, this paper aims to construct topology optimization for flow fields in VRFBs based on a three-dimensional model incorporating with electrochemical reaction kinetics. Referring to the models presented by several researchers \cite{shah2008dynamic,you2009simple,ma2011three}, a three-dimensional numerical model of a negative electrode in a VRFB is introduced and the electrolyte flow is assumed as stationary and isothermal Stokes flow for simplification. 
We demonstrate that the proposed approach enables the generation of a novel flow field configuration through the numerical example.
To confirm the performance of the topology optimized flow field, we investigate the mass transfer effect and overpotential of the topology optimized flow field in comparison with reference flow fields---parallel and interdigitated flow fields.
In addition, we discuss the power loss \cite{blanc2010understanding,xu2013numerical} in terms of polarization loss and pumping power at different operating conditions.

The reminder of this paper is organized as follows.
In Section 2, we introduce the mathematical model and assumptions of a VRFB.
In Section 3, we formulate a topology optimization problem that aims to maximize the mass transfer effect of a three-dimensional flow field in the VRFB and construct the optimization algorithm based on the use of mathematical programming and the finite element method (FEM).
In Section 4, we provide numerical examples and demonstrate the usefulness of the proposed approach.
Finally, Section 5 concludes this paper and summarizes the obtained results.

\section{Mathematical model}
\subsection{Model assumptions}
A schematic diagram of a typical redox flow battery is shown in Fig.~1. 
The positive and negative electrodes are separated by the ion exchange membrane, which only allows protons to penetrate. 
We suppose that the electrodes compose of the carbon fiber electrode and flow channel.
Note that the use of flow channel enables the reduction of pressure loss in comparison with the case of only using the carbon fiber electrode \cite{xu2013numerical}.
When the electrolyte stored in tanks circulates through the positive and negative electrode separately by pumps, the electric energy is released or stored by the electrochemical reactions in the electrodes. 
The main reactions can be described as follows:
\begin{align}
 & \text{Positive electrode: } \text{VO}^{2+} + \text{H}_{2}\text{O} \rightleftharpoons \text{VO}^{+}_{2} + 2\text{H}^{+} + \text{e}^{-} \\
 & \text{Negative electrode: } \text{V}^{3+} + \text{e} \rightleftharpoons \text{V}^{2+}
\end{align}

Furthermore, only the negative electrode is considered in this work and some basic assumptions are used for simplification as follows \cite{you2009simple}:
\begin{enumerate}[1.]
\item The electrolyte flow is treated as stationary, incompressible and isothermal Stokes flow.
\item The dilute-solution approximation is used in this numerical model.
\item The side reactions are neglected in electrochemical reactions.
\item The migration  phenomenon is ignored in the species transport process.
\end{enumerate}

\subsection{Governing equations}
Based on the assumptions presented in the previous section, governing equations incorporated in the numerical model are introduced here. The electrolyte flow passing through the flow channel in VRFBs can be described by Stokes equation and continuity equation as follows:
\begin{align}
   & -\nabla p +\mu\nabla^{2}\mathbf{u}= \mathbf{0}, \\
   & \nabla\cdot \mathbf{u} = 0,
\end{align} 
where $\mu$ is the viscosity of electrolyte flow, and $p(\mathbf{x})$ and $\mathbf{u}(\mathbf{x})$ are the pressure and velocity at position $\mathbf{x}$, respectively.

In addition to the flow channel, the porous electrode is also permeated with the electrolyte flow and the velocity in the electrode can be expressed by Darcy's law:
\begin{align}
  & \frac{\mu}{K}\mathbf{u} = -\nabla p,
\end{align}
where $K$ is the permeability coefficient, which can be described by the Kozeny-Carmen equation \cite{tomadakis2005viscous} as follows:
\begin{align}
  & K = \frac{d^{2}_\text{f}\epsilon^{3}}{16K_{\text{ck}}(1-\epsilon)^{2}},
\end{align}
where $d_{\text{f}}$ is the fiber diameter, $\epsilon$ is the porosity of electrodes,  $K_{\text{ck}}$ is the Carman-Kozeny constant described by the characteristic of the fibrous material.
In addition, instead of the electrolyte flow described by Stokes equation and Darcy's law respectively, 
Brinkman equation that combines Stokes equation with Darcy's law can be used to describe the electrolyte flow in the mixture of different porous medium generally and is formulated as
\begin{align}
   & -\nabla p + \mu\nabla^{2}\mathbf{u} +\mathbf{F}= \mathbf{0},
   \label{eq:br}
\end{align}
where $\mathbf{F}$ is the body force given by
\begin{align}
\mathbf{F} = - \alpha \mathbf{u},
\label{eq:f}
\end{align}
where $\alpha$ is the so-called inverse-permeability that is defined as $\alpha=\mu/K$ in the porous medium, while $\alpha=0$ in the pure fluid domain.

With the electrolyte flow containing vanadium species, the species transport needs to be considered in the numerical model, which can be expressed as follows:
\begin{align}
   & \mathbf{u}\cdot\nabla c_{i}-D^{\text{eff}}_{i}\nabla^{2}c_{i}=-s_{i},
\end{align} 
where $c_{i}$ is the concentration of vanadium species $i\in \{\text{V}^{2+}, \text{V}^{3+}$\}, and $s_{i}$ is the source term of species $i$ due to electrochemical reactions, $s_{\text{V}^{2+}}=j/F$ and $s_{\text{V}^{3+}}=-j/F$, where $j$ and $F$ are the transfer current density and the Faraday constant, respectively. 
In addition, $D^{\text{eff}}_{i}$ is the effective diffusion coefficient of species $i$ and is given by the Bruggemann correction, as follows:
\begin{align}
D^{\text{eff}}_{i} = \epsilon^{1.5}D_{i},
\end{align}
where $D_i$ is the diffusion coefficient of species $i$.
The charges in a VRFB conserve since the charge entering the electrolyte is balanced by the charge leaving the electrode, where the expression is shown as follows:
\begin{align}
   & \nabla\cdot \mathbf{i}_\text{e} + \nabla\cdot \mathbf{i}_\text{s} = 0,
   \label{eq:cc}
\end{align} 
where $\mathbf{i}_\text{e}$ is the ionic current density, and $\mathbf{i}_\text{s}$ is the electronic current density. However, when the charge flows from the electrolyte to the electrode, the electrochemical reactions should take place on the electrode surface, which must be also taken account into the charge conservation. 
With the electrochemical reaction kinetics, Eq.~(\ref{eq:cc}) can be rewritten as
\begin{align}
   & \nabla\cdot \mathbf{i}_\text{e} = -\nabla\cdot \mathbf{i}_\text{s} = j.
\end{align}
In addition, since the electrolyte is assumed to be electrically neutral and the migration term is ignored, the total ionic current density can be further expressed in terms of electric potential in the electrolyte as follows:
\begin{align}
   & \mathbf{i}_\text{e} = \sum_{i}\mathbf{i}_i = -\kappa^{\text{eff}}_{\text{e}}\nabla\phi_{\text{e}},
   \label{eq:ie}
   \\
   & \kappa^{\text{eff}}_{\text{e}} =  \frac{F^2}{RT}\sum_{i}z^{2}_{i}D^{\text{eff}}_{i}c_{i},
   \label{eq:ke}
\end{align} 
where  $\mathbf{i}_i$ is the ionic current density of species $i$, $\phi_\text{e}$ is the electric potential of the electrolyte,  $\kappa^{\text{eff}}_\text{e}$ is the effective conductivity of the electrolyte, $R$ is the gas constant, $T$ is the temperature, and $z_{i}$ is the valence.
The detailed derivation of Eqs.~(\ref{eq:ie}) and (\ref{eq:ke}) refers to the previous work by Shah et al.~\cite{shah2008dynamic}.

Likewise, the electronic current density can also be expressed in terms of electric potential in the electrode given by the Ohm's law, as follows:
\begin{align}
   & \mathbf{i}_\text{s} = -\sigma^{\text{eff}}_{\text{s}}\nabla\phi_{\text{s}}, \\
   & \sigma^{\text{eff}}_{\text{s}} = (1-\epsilon)^{1.5}\sigma_{\text{s}},
\end{align} 
where $\sigma_\text{s}$ is the conductivity of the solid material of the electrode, $\sigma^\text{eff}_\text{s}$ is the effective conductivity of the electrode, and $\phi_\text{s}$ is the electric potential of the electrode.
Note that the effective conductivity of the electrode is corrected by the Bruggemann correction.

The transfer current density $j$, which originates from electrochemical reactions, can be described using Butler-Volmer equation, as follows:
\begin{align}
&j = i_{0}\left[R^\text{sb}_{\text{V}^{3+}}\exp\left(-\frac{\alpha_\text{c}F\eta}{RT}\right)-R^\text{sb}_{\text{V}^{2+}}\exp\left(\frac{\alpha_\text{a}F\eta}{RT}\right)\right],\\
&i_{0} = aFk(c_{\text{V}^{2+}})^{\alpha_\text{c}}(c_{\text{V}^{3+}})^{\alpha_\text{a}},
\end{align}
where $i_{0}$ is the exchange current density, $\eta$ is the overpotential, $\alpha_\text{c}$ and $\alpha_\text{a}$ are the cathodic and anodic transfer coefficients, $k$ is the reaction rate constant, $a$ is the specific area of the electrode, and $R^\text{sb}_{i}=c^\text{s}_{i}/c_{i}$ is the ratio of the surface concentration of species $i$ to the bulk concentration in the negative electrode, in which $c^\text{s}_{i}$ is the species concentration at the surface of the negative electrode.

The overpotential in Butler-Volmer equation is the difference between the electrode potential and the electrolyte potential, which can be expressed as follows:
\begin{align}
&\eta = \phi_\text{s}-\phi_\text{e}-U,
\end{align}
where $U$ is the open-circuit potential in the negative electrode, which can be estimated by Nernst equation as follows:
\begin{align}
&U = U_{0}+\frac{RT}{F}\ln\left(\frac{c_{\text{V}^{3+}}}{c_{\text{V}^{2+}}}\right),
\end{align}
where $U_{0}$ is the equilibrium potential.

The species concentration at the electrode surface $c^\text{s}_{i}$ differs from the bulk concentration $c_{i}$ owing to the electrochemical reactions taking place at the electrode surface and the conductivity difference between the electrode and electrolyte. 
According to the previous research \cite{you2009simple}, $c^\text{s}_{i}$ are given by
\begin{align}
&c^\text{s}_{\text{V}^{2+}}=\frac{\overline{P}c_{\text{V}^{3+}}+(1+\overline{P})c_{\text{V}^{2+}}}{1+\overline{M}+\overline{P}},\\
&c^\text{s}_{\text{V}^{3+}}=\frac{\overline{M}c_{\text{V}^{2+}}+(1+\overline{M})c_{\text{V}^{3+}}}{1+\overline{M}+\overline{P}},
\label{eq:c2s}
\end{align}
where $\overline{M}$ and $\overline{P}$ are defined as follows:
\begin{align}
&\overline{M}=\frac{k}{k_\text{m}}(c_{\text{V}^{2+}})^{\alpha_\text{c}-1}(c_{\text{V}^{3+}})^{\alpha_\text{a}}\exp\left(\frac{\alpha_\text{a}F\eta}{RT}\right),\\
&\overline{P}=\frac{k}{k_\text{m}}(c_{\text{V}^{2+}})^{\alpha_\text{c}}(c_{\text{V}^{3+}})^{\alpha_\text{a}-1}\exp\left(-\frac{\alpha_\text{c}F\eta}{RT}\right),
\end{align}
where $k_\text{m}$ is the mass transfer coefficient, which can be estimated by the following equation \cite{schmal1986mass}:
\begin{align}
&k_\text{m} = 1.6 \times 10^{-4}|\mathbf{u}|^{0.4}.
\end{align}

\subsection{Boundary conditions}

\begin{figure}[t]
	\centering
\includegraphics[width=75mm]{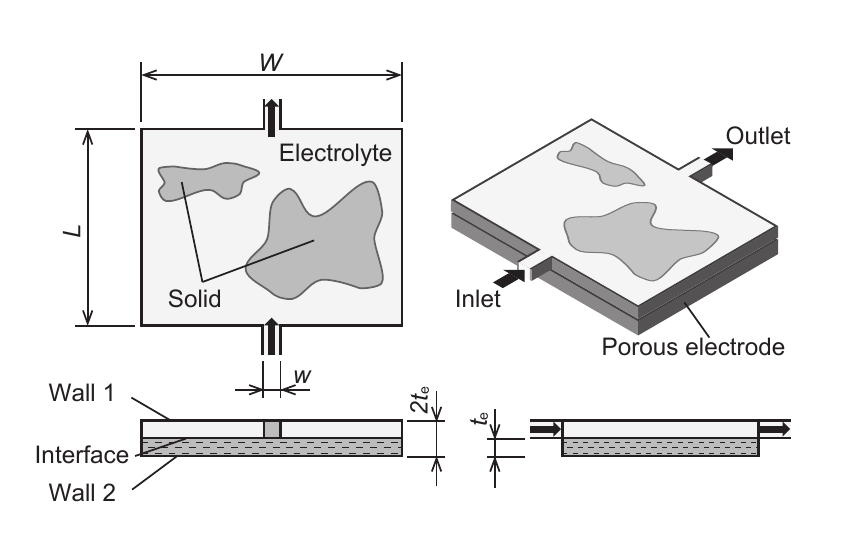}
	\vspace{-5mm}
	\caption{Schematic diagram of the analysis domain and boundary settings.}
	\label{fig: sec2_3_1}
\end{figure}

In this section, the boundary conditions used in the numerical model during topology optimization process are introduced. 
Figure \ref{fig: sec2_3_1} shows the schematic diagram of the analysis domain and boundary settings.
The pressure conditions are imposed on the inlet and the outlet, and the no-slip condition is applied on the remaining outer boundaries. The expressions are shown below:
\begin{alignat}{2}
&p = p_\text{in}  &&\ \ \text{on the inlet},\\
&p = p_\text{out} &&\ \ \text{on the outlet},\\
&\mathbf{u} =\mathbf{0} &&\ \ \text{on the remaining boundaries},
\end{alignat}
where $p_\text{in}$ is the given pressure value on the inlet, and $p_\text{out}$ is the given pressure value on the outlet.

\begin{table*}[t]
\centering  
\caption{Parameter settings of the electrode.}
	\vspace{-1mm}
\scalebox{0.8}{
\begin{tabular}[t]{lllll}
\hline
Parameter & Symbol & Value & Unit & Ref.\\
\hline
Porosity   & $\epsilon$   & 0.929 & - & \cite{you2009simple} \\
Specific surface area   & $a$   & $1.62\times10^{4}$ & m & \cite{you2009simple} \\
Carbon fiber diameter   & $d_\text{f}$   & $1.76\times10^{-5}$ & m & \cite{you2009simple} \\
Electronic conductivity of solid phase   & $\sigma_\text{s}$   & $1.0\times10^3$ & $\text{S m}^{-1}$ & \cite{you2009simple} \\
Kozeny-Carman constant   & $K_\text{ck}$   & $4.28$ & - & \cite{you2009simple} \\
Length   & $L$   & 0.1 & m & \cite{xu2013numerical} \\
Width   & $W$   & 0.1 & m & \cite{xu2013numerical} \\
Electrode thickness   & $t_\text{e}$   & $3.0\times10^{-3}$ & m & \cite{xu2013numerical} \\
\hline
\end{tabular}
}
\label{table:sec3_3_1}
\end{table*}
\begin{table*}[t]
\centering  
\caption{Parameter settings of the electrolyte.}
	\vspace{-1mm}
\scalebox{0.8}{
\begin{tabular}[t]{lllll}
\hline
Parameter & Symbol & Value & Unit & Ref.\\
\hline
Viscosity   & $\mu$   & $4.928\times10^{-3}$ & Pa s & \cite{you2009simple} \\
Initial vanadium $V^{2+}$concentration   & $c^\text{in}_{\text{V}^{2+}}$   & 750 & $\text{mol m}^{-3}$ & \cite{you2009simple} \\
Initial vanadium $V^{3+}$concentration   & $c^\text{in}_{\text{V}^{3+}}$   & 750 & $\text{mol m}^{-3}$ & \cite{you2009simple} \\
$V^{2+}$ diffusion coefficient   & $D_{\text{V}^{2+}}$   & $2.4\times10^{-4}$ & $\text{m}^{2}$ $\text{s}^{-1}$ & \cite{ma2011three} \\
$V^{3+}$ diffusion coefficient   & $D_{\text{V}^{3+}}$   & $2.4\times10^{-4}$ & $\text{m}^{2}$ $\text{s}^{-1}$ & \cite{ma2011three} \\
Ionic conductivity of electrolyte   & $\kappa_\text{e}$   & 7.8 & $\text{S m}^{-1}$ & Estimated \\

\hline
\end{tabular}
}
\label{table:sec3_3_2}
\end{table*}

For species conservation, the given concentration of each species is applied on the inlet and the diffusive fluxes of each species are set to zero on the outlet. Besides, the no-flux condition is also applied on the remaining boundaries. The boundary conditions for species conservation are shown below:
\begin{alignat}{2}
&c_{i} = c^\text{in}_{i}  &&\ \ \text{on the inlet},\\
&-D^\text{eff}_{i}\nabla c_{i}\cdot\mathbf{n}=0 &&\ \ \text{on the remaining boundaries}.
\label{eq:nof}
\end{alignat}
For charge conservation, the VRFB is assumed to be operated in galvanostatic situation. During charging process, the flux conditions in the negative electrode can be described as follows: 
\begin{alignat}{2}
&-\sigma^\text{eff}_\text{s}\nabla \phi_\text{s}\cdot \mathbf{n} = -{\frac{I}{A}} &&\ \ \text{on the wall 1},\\
&-\kappa^\text{eff}_\text{e}\nabla \phi_\text{e}\cdot \mathbf{n} = {\frac{I}{A}} &&\ \ \text{on the wall 2},\\
& \phi_\text{s} = 0 &&\ \ \text{on the interface},
\end{alignat}
where $I$ is the applied current, and $A$ is the electrode surface area.
Note that the remaining boundary conditions for $\phi_\text{s}$ and $\phi_\text{e}$ are the no-flux conditions as with Eq.~(\ref{eq:nof}).

\section{Topology optimization for flow fields}

\subsection{Concepts of topology optimization for fluid problems}
\begin{figure*}[t]
	\centering
	\includegraphics[width=150mm]{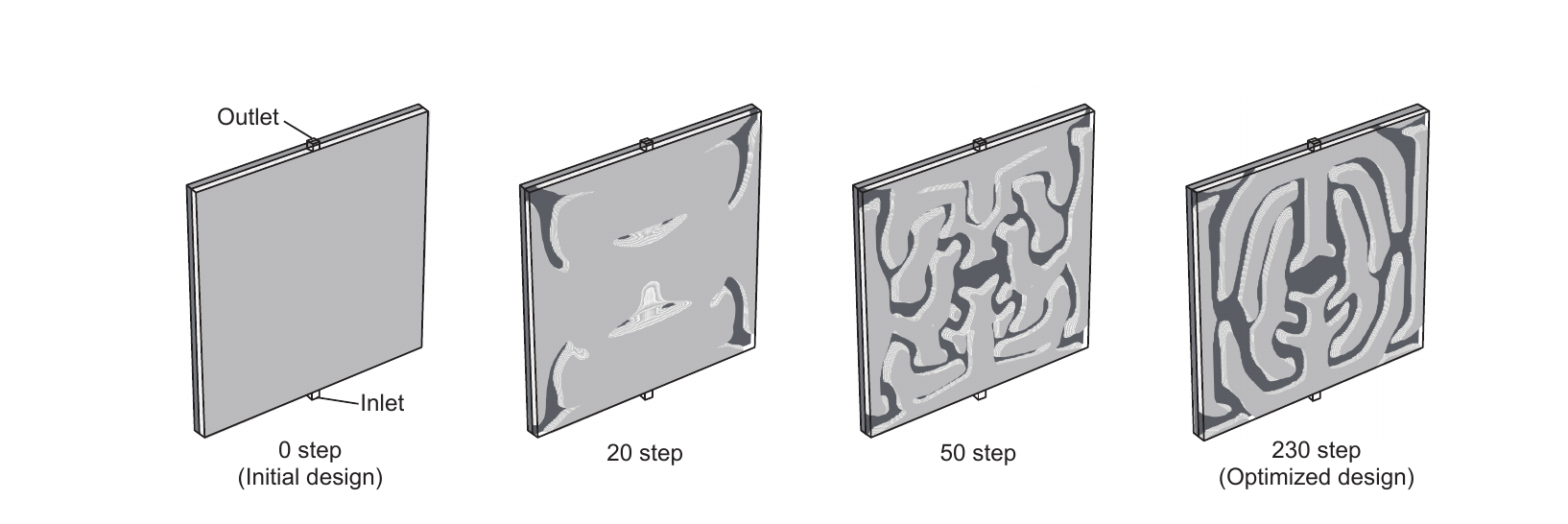}
	\vspace{-5mm}
\caption{Iteration history of topology optimized flow field expressed as $\rho\ge 0.5$.}
	\label{fig: sec4_1_1}
\end{figure*}

\begin{table*}[t]
\centering  
\caption{Parameter settings of the electrochemical reaction model.}
	\vspace{-1mm}
\scalebox{0.8}{
\begin{tabular}[t]{lllll}
\hline
Parameter & Symbol & Value & Unit & Ref.\\
\hline
Standard reaction rate  & $k_\text{c}$   & $1.7\times10^{-7}$ & $\text{m s}^{-1}$ & \cite{you2009simple} \\
Cathodic transfer coefficient   & $\alpha_\text{c}$   & 0.5 & - & Assumed \\
Anodic transfer coefficient   & $\alpha_\text{a}$   & 0.5 &- & Assumed \\
Equilibrium    & $U_{0}$   & -0.255 & V & \cite{ma2011three} \\

\hline
\end{tabular}
}
\label{table:sec3_3_3}
\end{table*}

Topology optimization aims to obtain the improved structural design in a specific domain with a given objective function and constraints. The main idea is to formulate a topology optimization problem as a material distribution problem, where the expression for material distribution in a fixed design domain $D$ is defined as follows \cite{bendsoe2003topology}: 
\begin{eqnarray}
\chi(\mathbf{x})=
\begin{cases}
1 &\text{if }\mathbf{x}\in \Omega,\cr 0 &\text{if }\mathbf{x}\in  D\backslash\Omega, \end{cases}
\end{eqnarray}
where $\mathbf{x}$ is the position in $D$, and $\Omega$ is the design domain in $D$. In this expression, $\chi(\mathbf{x}) = 1$ and $\chi(\mathbf{x}) = 0$ represent the material and void at $\mathbf{x}$, respectively. 
Since the characteristic function $\chi$ is a discontinuous function, topology optimization problems typically require relaxation techniques for numerical treatment.
As the popular and simple way for relaxing topology optimization problems, the density approach \cite{bendsoe1989optimal} replaces the characteristic function with a continuous function, $0 \le \rho(\mathbf{\mathbf{x}}) \le 1$, which is also used in this paper.

To determine which points in $D$ should be fluid or solid, based on the previous research dealing with fluid topology optimization \cite{borrvall2003topology}, the body force in Eq.~(\ref{eq:f}) is redefined using the fictitious body force, $\mathbf{F}^{\text{fic}}$, as follows:
\begin{align}
   &\mathbf{F}^{\text{fic}}= -\alpha^{\text{fic}}_{\rho} \mathbf{u} \quad\text{with }\ \alpha^{\text{fic}}_{\rho} = \frac{q(1-\rho)}{\rho+q}\alpha^{\text{fic}},
   \label{eq:f_fic}
\end{align}
where $\alpha^{\text{fic}}$ is the fictitious inverse-permeability used for expressing the solid domain $D\setminus\Omega$ as with the previous work \cite{borrvall2003topology}, and $q$ is a tuning parameter for controlling the convexity of $\alpha^{\text{fic}}_{\rho}$. 
In addition, $\rho = 1$ represents the fluid domain with $\alpha^{\text{fic}}_{\rho}=0$, and $\rho = 0$ represents the solid domain with $\alpha^{\text{fic}}_{\rho}=\alpha^{\text{fic}}\gg 1$. In the solid domain, since the fictitious body force is large enough compared to the fluid domain, it is difficult for fluid to pass through the solid domain. Therefore, according to the sensitivity information of the objective function, the fictitious body force is determined at each point in $D$ so that the structural design of flow channel can be obtained. 
Note that the fictitious inverse permeability, $\alpha^{\text{fic}}$, is different from $\alpha$ in Eq.~(\ref{eq:br}).
That is, the former is used for expressing the solid domain in the fixed design domain $D$, whereas the latter is used for expressing the porous electrode that is the non-design domain. 
In this study, $q$ and $\alpha^\text{fic}$ are set to $0.01$ and $5\alpha$, respectively.

\subsection{Description of optimization problem}

In VRFBs, the mobility of ions in electrolyte is poor compared to the electrons, which means it is more difficult to reach the electrode surface for ions. 
With the flow channel embedded in a VRFB, the mass transfer effect is improved so that the concentration of reactants at the electrode surface will also increase. Therefore, whether the mass transfer effect is improved can be estimated by the concentration of the reactants at the electrode surface. During charging process, the optimization problem can be defined as a maximization problem of average concentration of oxidized reactants at the electrode surface in a negative electrode. The expressions of this optimization problem are shown as follows:
\begin{align}
&\begin{array}{ll}
\displaystyle\underset{\rho}{\text{maximize }}\ F = \int_{D}c^\text{s}_{\text{V}^{3+}}\text{d}\Omega\bigg/\int_{D}\text{d}\Omega,\\
\text{subject to }\ 0 \le \rho(\mathbf{x}) \le 1\quad\text{for }\ \forall\mathbf{x}\in D.
\end{array}
\label{eq:optp}
\end{align}
Note that the detailed expressions of $c^\text{s}_{\text{V}^{3+}}$ can be found in Eq.~(\ref{eq:c2s}).

\subsection{Numerical implementation}

The governing equations are solved by using the package COMSOL Multiphysics$^{\textregistered}$, which is based on the FEM. The optimization algorithm is constructed on the basis of mathematical programming and is briefly enumerated as follows:
\begin{description}
\item[{\it Step 1.}]The design variables and all of the parameters shown in Table \ref{table:sec3_3_1}--\ref{table:sec3_3_4} are initialized.
\item[{\it Step 2.}]The objective function $F$ in (\ref{eq:optp}) is evaluated by solving the governing equation via the FEM.
\item[{\it Step 3.}] If the objective function is converged, the iteration will terminate. Otherwise, the sensitivities---gradient of objective function with respect to the design variables---are calculated.
\item[{\it Step 4.}] The design variables are redistributed in the fixed design domain $D$ using sequential linear programming (SLP), and the iteration will return to the second step.
\end{description}

We utilize a partial differential equation (PDE)-based filter \cite{kawamoto2011heaviside} for ensuring the smoothness of the design variables \cite{yaji2018topology}. 
In addition, we use the adjoint method that enables the derivation of the sensitivities without depending on the number of design variables. 
The detailed concepts and formulation of the adjoint method can be seen in the literatures on structural optimization \cite{haftka1992elements,bendsoe2003topology}.

\begin{figure*}[t]
	\centering
\includegraphics[width=150mm]{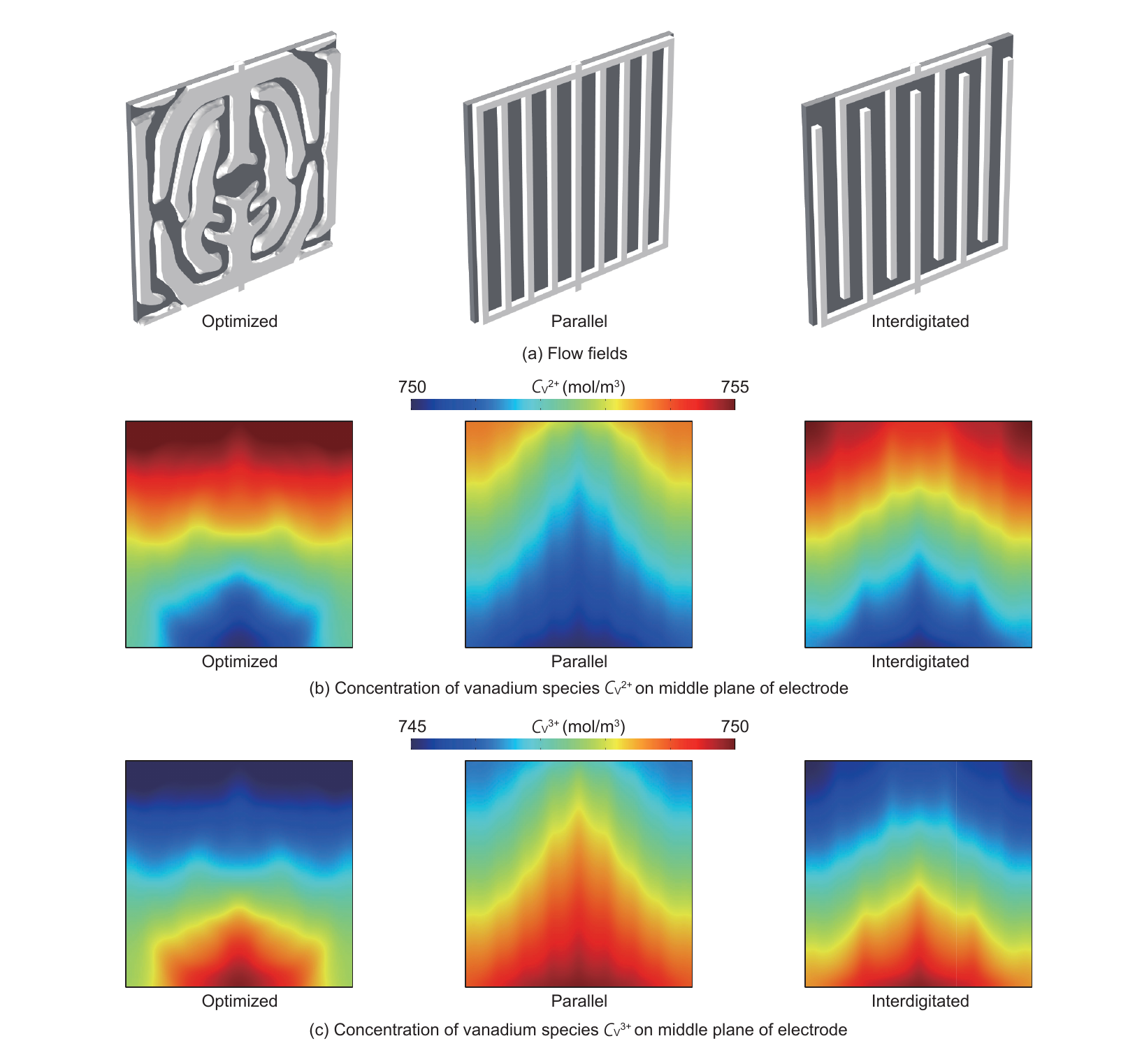}
	\vspace{-5mm}
	\caption{Comparison results for different flow fields: (a) Flow field configurations; (b) Distributions of vanadium species $c_{\text{V}^{2+}}$ on middle plane of electrode; (c) Distributions of vanadium species $c_{\text{V}^{3+}}$ on middle plane of electrode.}
	\label{fig: sec4_2_2}
\end{figure*}
\begin{figure*}[t]
	\centering
\includegraphics[width=150mm]{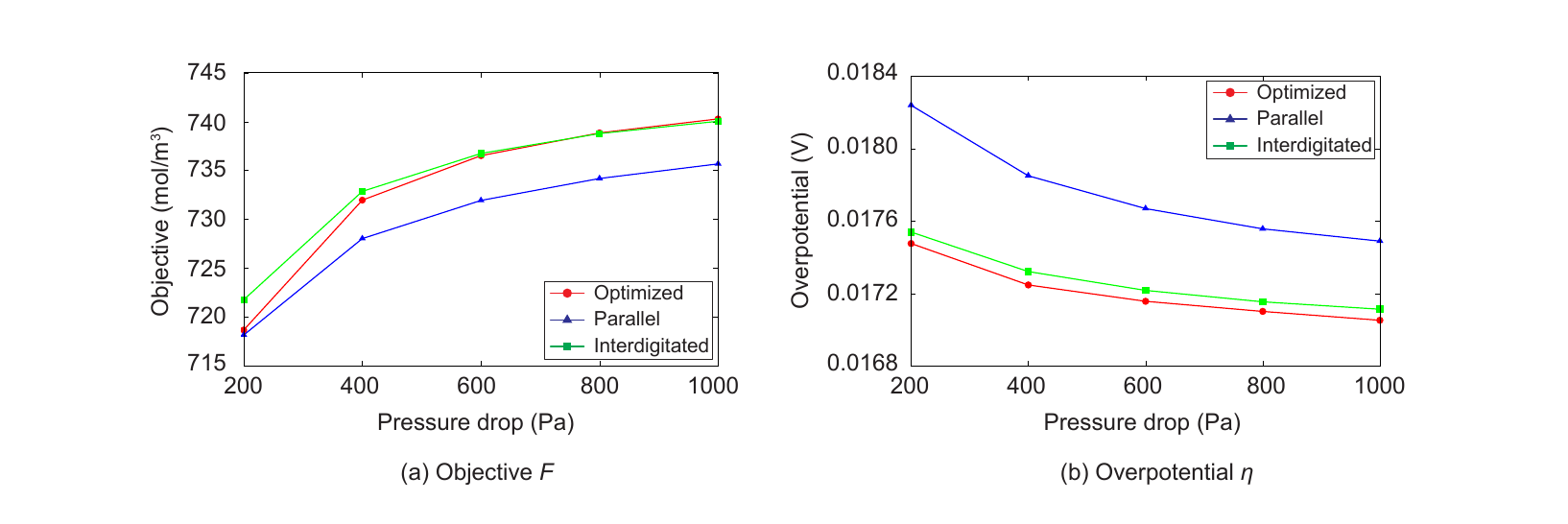}
	\vspace{-5mm}
	\caption{Values of objective and overpotential at variant pressure drop for different flow fields.}
	\label{fig: sec4_2_1}
\end{figure*}
\begin{figure*}[t]
	\centering
\includegraphics[width=150mm]{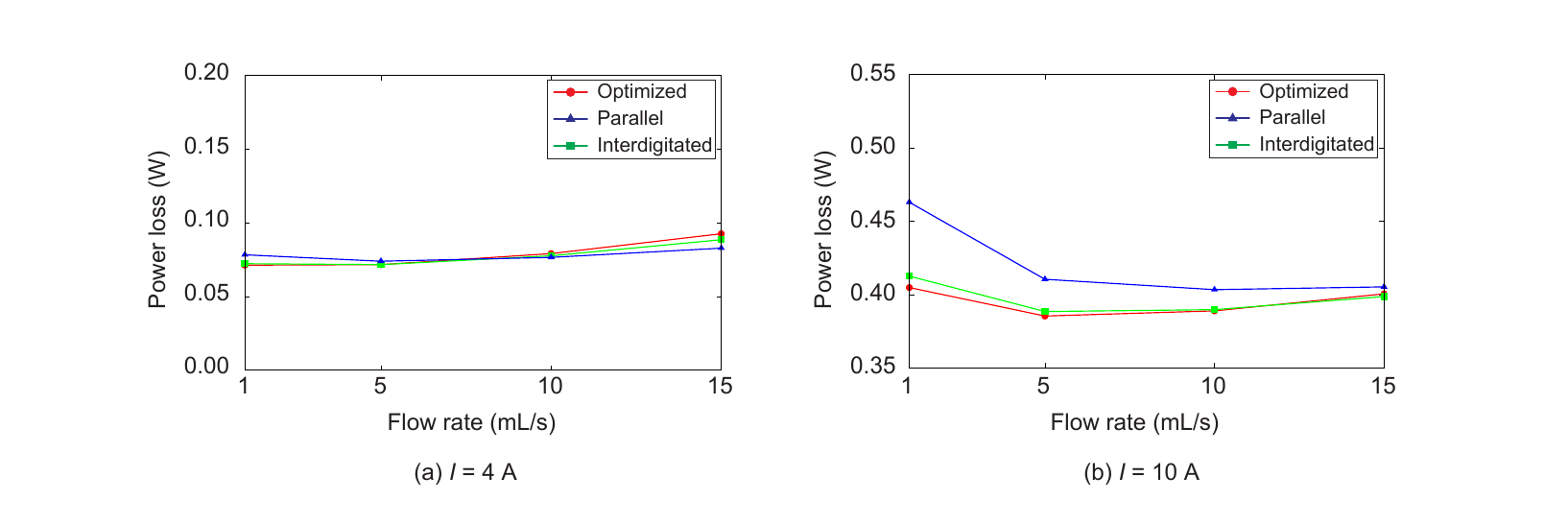}
	\vspace{-5mm}
	\caption{Power loss with different flow fields at $\epsilon = 0.929$.}
	\label{fig: sec4_3_1}
\end{figure*}
\begin{figure}[t]
	\centering
\includegraphics[width=75mm]{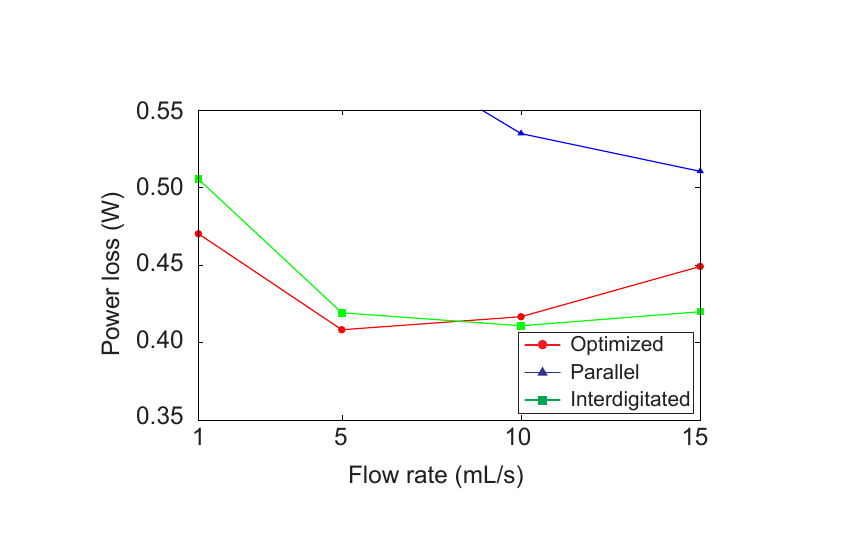}
	\vspace{-5mm}
	\caption{Power loss with different flow fields at $\epsilon = 0.68$ and $I = 10$~A.}
	\label{fig: sec4_3_3}
\end{figure}

\section{Numerical examples}

\subsection{Topology optimized flow field}
\begin{table}[b]
\centering  
\caption{Operating parameter settings.}
	\vspace{-1mm}
\scalebox{0.8}{
\begin{tabular}[t]{lllll}
\hline
Parameter & Symbol & Value & Unit & Ref.\\
\hline
Temperature  & $T$   & $298$ & K & \cite{you2009simple} \\
Inlet pressure   & $P_\text{in}$   & $1.0\times10^3$ & Pa & Assumed \\
Outlet pressure   & $P_\text{out}$   & 0 & Pa & Assumed \\
Applied current    & $I$   & 4.0 & A & \cite{xu2013numerical} \\

\hline
\end{tabular}
}
\label{table:sec3_3_4}
\end{table}
\begin{table}[b]
\centering  
\caption{Pressure drop (Pa) in the different flow fields.}
	\vspace{-1mm}
\scalebox{0.8}{
\begin{tabular}[t]{lllll}
\hline
                         & 1 mL/s    & 5 mL/s    & 10 mL/s     & 15 mL/s\\
\hline
Parallel      & 51   & 262 & 541   & 833 \\
Interdigitated       & 83   & 427 & 876   & 1343 \\
Optimized   & 105  & 531 & 1075 &1630 \\
\hline
\end{tabular}
}
\label{table:sec4_3_1}
\end{table}

Figure \ref{fig: sec4_1_1} shows the iteration history of topology optimized flow field, in which the electrolyte domain is expressed as the isosurface of $\rho\ge0.5$.
The analysis domain is discretized using $2.4\times 10^5$ hexahedral elements for all variables in this study.
In the topology optimized design shown in Fig.~\ref{fig: sec4_1_1}, as with the interdigitated flow field, some of the flow channels are not connected in the topology optimized flow field, where the disconnectivity of the flow field can enhance the mass transfer effect so that the mass transfer loss can be further reduced. 

\subsection{Effect of flow field design}

To compare the differences between the concentration of vanadium species with different flow fields, Fig.~\ref{fig: sec4_2_2} shows the geometric model of the flow fields and the distributions of $c_{\text{V}^{2+}}$ and $c_{\text{V}^{3+}}$ on the plane located in the middle of the electrode.
The geometric parameters of the electrode refers to Table \ref{table:sec3_3_1}. 
Besides, the width and thickness of the flow channel is 3~mm and the interval between the branches of the flow field is 9~mm. 
As shown in Fig.~\ref{fig: sec4_2_2}, the under-rib convection of the parallel flow field is weak, especially on both sides of the electrode. The velocity of electrolyte flow to both sides is small in the electrode since it is difficult for the parallel flow field to distribute the electrolyte to both sides of the electrode effectively. In contrast,  the under-rib convection is strong in the interdigitated flow field and topology optimized flow field. 

Figure \ref{fig: sec4_2_1} shows that the value of objective function in the optimization problem described as (\ref{eq:optp}), and overpotential in the topology optimized flow field and interdigitated flow filed indicate higher performance than those of the parallel flow field at constant pressure drop. 
Since the under-rib convection is weak in the parallel flow field, the value of objective function with the parallel flow field is the lowest. 
Besides, the electrode surface concentration decreases, as the pressure drop decreases for all the flow fields.
In other words, the species will be more difficult to reach the electrode surface at low flow rate, which means the mass transfer effect is dominated by the flow channel.
Accordingly, the overpotential with topology optimized flow field is the lowest among the flow fields due to the strongest mass transfer effect caused by the topology optimized flow field.
However, it should be noted that the overall performances of the optimized flow field and the interdigitated flow field are almost same in this numerical example.

\subsection{Power loss}

In a VRFB system, the polarization loss is caused by the need for extra voltage to occur electrochemical reactions. 
In addition to the polarization loss, the pumping power should be also taken into account for the evaluation of VRFB system.
We therefore introduce the following evaluation index:
\begin{align}
P_{\text{loss}} = I\eta + Q\Delta P,
\label{eq:loss}
\end{align}
where $P_{\text{loss}}$ is the sum of the polarization loss ($I\eta$), and the pumping power loss ($Q\Delta P$), $Q$ is the flow rate, and $\Delta P$ is the pressure drop.
Note that, for brevity in this study, the evaluation index in Eq.~(\ref{eq:loss}) is defined as a simple expression using the dominant factors, whereas the performance of VRFB systems relates with various factors on the authority of the previous works \cite{blanc2010understanding,xu2013numerical}.

Figure \ref{fig: sec4_3_1} shows the power loss corresponding to variant flow rates. At low flow rate of 1 mL$/$s, the power loss with the topology optimized flow field is the lowest one at current of 4 A or 10 A. Although the pressure drop in the topology optimized flow field shown in Table \ref{table:sec4_3_1} is the highest, the pumping power is trivial compared to the polarization loss, which results in the lowest power loss with the topology optimized flow field. At high flow rate of 15~mL$/$s, the parallel flow field shows the lowest power loss at current of 4~A. However, the power loss with these flow fields is almost the same at current of 10~A. At high applied current, the electrochemical reactions at the electrode surface become strong and the polarization loss will increase. Therefore, with the strongest mass transfer effect caused by the topology optimized flow field, it is obvious that the difference of the total power loss between these flow fields decreases as the applied current density increases at high flow rates, even if the pumping power increases as the flow rate increases. In contrast, at flow rate of 1~mL$/$s, the differences in the power loss between the flow fields increases as the applied current increases because the power loss is mainly from the polarization loss at current of 4~A and the polarization loss will be more dominant in power loss as the current increases from 4~A to 10~A. The results demonstrate the topology optimized flow field is more useful for VRFBs at high applied current.

Figure \ref{fig: sec4_3_3} shows the comparison results of the power loss at different operating condition at $\epsilon=0.68$ and $I=10$~A. For all the flow fields, the power loss increases in comparison with the case of $\epsilon=0.929$ in Fig.~\ref{fig: sec4_3_1}. 
At flow rate of 1~mL$/$s, with the weaker mass transfer effect, the power loss of the interdigitated flow field and the topology optimized flow field increases from 0.413 to 0.505 and from 0.405 to 0.470 as the porosity decreases from 0.929 to 0.68. 
The differences between the flow fields become larger at low porosity. At flow rate of 15~mL$/$s, since the pressure drop in the topology optimized flow field is larger than that of the interdigitated flow field and the polarization loss falls slowly when the flow rate increases from 1 mL$/$s to 15 mL$/$s, the total power loss with interdigitated flow field is less than the loss with the topology optimized flow field. For low porosity, the polarization loss and pumping power will be more dominant at low and high flow rate respectively. 

\section{Conclusion}
In this paper, we proposed a topology optimization method for the flow field in a VRFB based on a three-dimensional numerical model. 
We demonstrated that the optimization problem can be formulated as a maximization problem of the electrode surface concentration of oxidized reactants during the charging process.
Based on the proposed formulations, we derived a novel flow field configuration and evaluated its performances in comparison with reference flow fields---parallel and interdigitated flow fields.
As a result, we verified that a VRFB with the topology optimized flow field can achieve almost same performances with those of the interdigitated flow field, which has stronger mass transfer effect than the parallel flow field.

%
We further investigated the power loss with the parallel, interdigitated and topology optimized flow fields under different operating conditions.
The results demonstrated that a VRFB with the topology optimized flow field is more suitable for a VRFB at high applied current density.
We confirmed that the topology optimized flow field has a potential to be an alternative design candidate when dealing with low porosity electrode, as the difference of the power loss with respect to the flow field configurations is notably observed in comparison with the case of high porosity electrode. 

\section*{Acknowledgments}
This work is partially supported by a research grant from The Mazda Foundation.
\color{black}

\bibliography{mybibfile}

\begin{thebibliography}{10}
\expandafter\ifx\csname url\endcsname\relax
  \def\url#1{\texttt{#1}}\fi
\expandafter\ifx\csname urlprefix\endcsname\relax\def\urlprefix{URL }\fi
\expandafter\ifx\csname href\endcsname\relax
  \def\href#1#2{#2} \def\path#1{#1}\fi

\bibitem{alotto2014redox}
P.~Alotto, M.~Guarnieri, F.~Moro, {Redox flow batteries for the storage of
  renewable energy: A review}, Renewable and Sustainable Energy Reviews 29
  (2014) 325--335.

\bibitem{milshtein2017quantifying}
J.~D. Milshtein, K.~M. Tenny, J.~L. Barton, J.~Drake, R.~M. Darling, F.~R.
  Brushett, Quantifying mass transfer rates in redox flow batteries, Journal of
  The Electrochemical Society 164~(11) (2017) E3265--E3275.

\bibitem{aaron2012dramatic}
D.~Aaron, Q.~Liu, Z.~Tang, G.~Grim, A.~Papandrew, A.~Turhan, T.~Zawodzinski,
  M.~Mench, Dramatic performance gains in vanadium redox flow batteries through
  modified cell architecture, Journal of Power Sources 206 (2012) 450--453.

\bibitem{sun1992modification}
B.~Sun, M.~Skyllas-Kazacos, {Modification of graphite electrode materials for
  vanadium redox flow battery application---I. Thermal treatment},
  Electrochimica Acta 37~(7) (1992) 1253--1260.

\bibitem{wang2007investigation}
W.~Wang, X.~Wang, {Investigation of Ir-modified carbon felt as the positive
  electrode of an all-vanadium redox flow battery}, Electrochimica Acta 52~(24)
  (2007) 6755--6762.

\bibitem{chen2013optimizing}
D.~Chen, M.~A. Hickner, E.~Agar, E.~C. Kumbur, Optimizing membrane thickness
  for vanadium redox flow batteries, Journal of Membrane Science 437 (2013)
  108--113.

\bibitem{li2013bismuth}
B.~Li, M.~Gu, Z.~Nie, Y.~Shao, Q.~Luo, X.~Wei, X.~Li, J.~Xiao, C.~Wang,
  V.~Sprenkle, W.~Wang, Bismuth nanoparticle decorating graphite felt as a
  high-performance electrode for an all-vanadium redox flow battery, Nano
  Letters 13~(3) (2013) 1330--1335.

\bibitem{li2011graphite}
W.~Li, J.~Liu, C.~Yan, Graphite--graphite oxide composite electrode for
  vanadium redox flow battery, Electrochimica Acta 56~(14) (2011) 5290--5294.

\bibitem{zhou2017critical}
X.~Zhou, T.~Zhao, L.~An, Y.~Zeng, L.~Wei, Critical transport issues for
  improving the performance of aqueous redox flow batteries, Journal of Power
  Sources 339 (2017) 1--12.

\bibitem{xu2013numerical}
Q.~Xu, T.~Zhao, P.~Leung, Numerical investigations of flow field designs for
  vanadium redox flow batteries, Applied Energy 105 (2013) 47--56.

\bibitem{tsushima2014efficient}
S.~Tsushima, F.~Kondo, S.~Sasaki, S.~Hirai, Efficient utilization of the
  electrodes in a redox flow battery by modifying flow field and electrode
  morphology, Proceedings of 15th International Heat Transfer Conference (2014)
  doi:10.1615/IHTC15.ecs.009326.

\bibitem{darling2014influence}
R.~M. Darling, M.~L. Perry, The influence of electrode and channel
  configurations on flow battery performance, Journal of The Electrochemical
  Society 161~(9) (2014) A1381--A1387.

\bibitem{houser2016influence}
J.~Houser, J.~Clement, A.~Pezeshki, M.~M. Mench, Influence of architecture and
  material properties on vanadium redox flow battery performance, Journal of
  Power Sources 302 (2016) 369--377.

\bibitem{bendsoe1988generating}
M.~P. Bends{\o}e, N.~Kikuchi, Generating optimal topologies in structural
  design using a homogenization method, Computer Methods in Applied Mechanics
  and Engineering 71~(2) (1988) 197--224.

\bibitem{bendsoe2003topology}
M.~P. Bends{\o}e, O.~Sigmund, Topology optimization: theory, methods, and
  applications, Springer, 2003.

\bibitem{bendsoe1989optimal}
M.~P. Bends{\o}e, Optimal shape design as a material distribution problem,
  Structural Optimization 1~(4) (1989) 193--202.

\bibitem{diaz1992solutions}
A.~R. Diaz, N.~Kikuchi, Solutions to shape and topology eigenvalue optimization
  problems using a homogenization method, International Journal for Numerical
  Methods in Engineering 35~(7) (1992) 1487--1502.

\bibitem{ma1995topological}
Z.-D. Ma, N.~Kikuchi, H.-C. Cheng, Topological design for vibrating structures,
  Computer Methods in Applied Mechanics and Engineering 121~(1-4) (1995)
  259--280.

\bibitem{li1999shape}
Q.~Li, G.~P. Steven, O.~M. Querin, Y.~Xie, Shape and topology design for heat
  conduction by evolutionary structural optimization, International Journal of
  Heat and Mass Transfer 42~(17) (1999) 3361--3371.

\bibitem{iga2009topology}
A.~Iga, S.~Nishiwaki, K.~Izui, M.~Yoshimura, Topology optimization for thermal
  conductors considering design-dependent effects, including heat conduction
  and convection, International Journal of Heat and Mass Transfer 52~(11-12)
  (2009) 2721--2732.

\bibitem{nomura2007structural}
T.~Nomura, K.~Sato, K.~Taguchi, T.~Kashiwa, S.~Nishiwaki, Structural topology
  optimization for the design of broadband dielectric resonator antennas using
  the finite difference time domain technique, International Journal for
  Numerical Methods in Engineering 71~(11) (2007) 1261--1296.

\bibitem{yamasaki2011level}
S.~Yamasaki, T.~Nomura, A.~Kawamoto, K.~Sato, S.~Nishiwaki, A level set-based
  topology optimization method targeting metallic waveguide design problems,
  International Journal for Numerical Methods in Engineering 87~(9) (2011)
  844--868.

\bibitem{sigmund2001design}
O.~Sigmund, {Design of multiphysics actuators using topology
  optimization---Part I: One-material structures}, Computer Methods in Applied
  Mechanics and Engineering 190~(49-50) (2001) 6577--6604.

\bibitem{iwai2011power}
H.~Iwai, A.~Kuroyanagi, M.~Saito, A.~Konno, H.~Yoshida, T.~Yamada,
  S.~Nishiwaki, Power generation enhancement of solid oxide fuel cell by
  cathode--electrolyte interface modification in mesoscale assisted by level
  set-based optimization calculation, Journal of Power Sources 196~(7) (2011)
  3485--3495.

\bibitem{song20132d}
X.~Song, A.~Diaz, A.~Benard, J.~Nicholas, A 2d model for shape optimization of
  solid oxide fuel cell cathodes, Structural and Multidisciplinary Optimization
  47~(3) (2013) 453--464.

\bibitem{borrvall2003topology}
T.~Borrvall, J.~Petersson, {Topology optimization of fluids in Stokes flow},
  International Journal for Numerical Methods in Fluids 41~(1) (2003) 77--107.

\bibitem{gersborg2005topology}
A.~Gersborg-Hansen, O.~Sigmund, R.~B. Haber, Topology optimization of channel
  flow problems, Structural and Multidisciplinary Optimization 30~(3) (2005)
  181--192.

\bibitem{olesen2006high}
L.~H. Olesen, F.~Okkels, H.~Bruus, {A high-level programming-language
  implementation of topology optimization applied to steady-state
  Navier--Stokes flow}, International Journal for Numerical Methods in
  Engineering 65~(7) (2006) 975--1001.

\bibitem{kubo2017level}
S.~Kubo, K.~Yaji, T.~Yamada, K.~Izui, S.~Nishiwaki, A level set-based topology
  optimization method for optimal manifold designs with flow uniformity in
  plate-type microchannel reactors, Structural and Multidisciplinary
  Optimization 55~(4) (2017) 1311--1327.

\bibitem{yoon2016topology}
G.~H. Yoon, {Topology optimization for turbulent flow with Spalart--Allmaras
  model}, Computer Methods in Applied Mechanics and Engineering 303 (2016)
  288--311.

\bibitem{dilgen2018topology}
C.~B. Dilgen, S.~B. Dilgen, D.~R. Fuhrman, O.~Sigmund, B.~S. Lazarov, Topology
  optimization of turbulent flows, Computer Methods in Applied Mechanics and
  Engineering 331 (2018) 363--393.

\bibitem{yoon2010topology}
G.~H. Yoon, Topology optimization for stationary fluid--structure interaction
  problems using a new monolithic formulation, International Journal for
  Numerical Methods in Engineering 82~(5) (2010) 591--616.

\bibitem{jenkins2015level}
N.~Jenkins, K.~Maute, Level set topology optimization of stationary
  fluid-structure interaction problems, Structural and Multidisciplinary
  Optimization 52~(1) (2015) 179--195.

\bibitem{matsumori2013topology}
T.~Matsumori, T.~Kondoh, A.~Kawamoto, T.~Nomura, Topology optimization for
  fluid--thermal interaction problems under constant input power, Structural
  and Multidisciplinary Optimization 47~(4) (2013) 571--581.

\bibitem{yaji2015topology}
K.~Yaji, T.~Yamada, S.~Kubo, K.~Izui, S.~Nishiwaki, A topology optimization
  method for a coupled thermal--fluid problem using level set boundary
  expressions, International Journal of Heat and Mass Transfer 81 (2015)
  878--888.

\bibitem{yaji2018large}
K.~Yaji, M.~Ogino, C.~Chen, K.~Fujita, Large-scale topology optimization
  incorporating local-in-time adjoint-based method for unsteady thermal-fluid
  problem, Structural and Multidisciplinary Optimization (2018)
  doi:10.1007/s00158--018--1922--6.

\bibitem{alexandersen2014topology}
J.~Alexandersen, N.~Aage, C.~S. Andreasen, O.~Sigmund, Topology optimisation
  for natural convection problems, International Journal for Numerical Methods
  in Fluids 76~(10) (2014) 699--721.

\bibitem{coffin2016level}
P.~Coffin, K.~Maute, A level-set method for steady-state and transient natural
  convection problems, Structural and Multidisciplinary Optimization 53~(5)
  (2016) 1047--1067.

\bibitem{alexandersen2016large}
J.~Alexandersen, O.~Sigmund, N.~Aage, Large scale three-dimensional topology
  optimisation of heat sinks cooled by natural convection, International
  Journal of Heat and Mass Transfer 100 (2016) 876--891.

\bibitem{kontoleontos2013adjoint}
E.~Kontoleontos, E.~Papoutsis-Kiachagias, A.~Zymaris, D.~Papadimitriou,
  K.~Giannakoglou, Adjoint-based constrained topology optimization for viscous
  flows, including heat transfer, Engineering Optimization 45~(8) (2013)
  941--961.

\bibitem{dilgen2018density}
S.~B. Dilgen, C.~B. Dilgen, D.~R. Fuhrman, O.~Sigmund, B.~S. Lazarov, Density
  based topology optimization of turbulent flow heat transfer systems,
  Structural and Multidisciplinary Optimization (2018)
  doi:10.1007/s00158--018--1967--6.

\bibitem{yaji2018topology}
K.~Yaji, S.~Yamasaki, S.~Tsushima, T.~Suzuki, K.~Fujita, Topology optimization
  for the design of flow fields in a redox flow battery, Structural and
  Multidisciplinary Optimization 57~(2) (2018) 535--546.

\bibitem{shah2008dynamic}
A.~Shah, M.~Watt-Smith, F.~Walsh, A dynamic performance model for redox-flow
  batteries involving soluble species, Electrochimica Acta 53~(27) (2008)
  8087--8100.

\bibitem{you2009simple}
D.~You, H.~Zhang, J.~Chen, A simple model for the vanadium redox battery,
  Electrochimica Acta 54~(27) (2009) 6827--6836.

\bibitem{ma2011three}
X.~Ma, H.~Zhang, F.~Xing, A three-dimensional model for negative half cell of
  the vanadium redox flow battery, Electrochimica Acta 58 (2011) 238--246.

\bibitem{blanc2010understanding}
C.~Blanc, A.~Rufer, Understanding the vanadium redox flow batteries, in: Paths
  to Sustainable Energy, InTech, 2010.

\bibitem{tomadakis2005viscous}
M.~M. Tomadakis, T.~J. Robertson, Viscous permeability of random fiber
  structures: comparison of electrical and diffusional estimates with
  experimental and analytical results, Journal of Composite Materials 39~(2)
  (2005) 163--188.

\bibitem{schmal1986mass}
D.~Schmal, J.~Van~Erkel, P.~Van~Duin, Mass transfer at carbon fibre electrodes,
  Journal of Applied Electrochemistry 16~(3) (1986) 422--430.

\bibitem{kawamoto2011heaviside}
A.~Kawamoto, T.~Matsumori, S.~Yamasaki, T.~Nomura, T.~Kondoh, S.~Nishiwaki,
  {Heaviside projection based topology optimization by a PDE-filtered scalar
  function}, Structural and Multidisciplinary Optimization 44~(1) (2011)
  19--24.

\bibitem{haftka1992elements}
R.~T. Haftka, Z.~G{\"u}rdal, Elements of structural optimization, 3rd edn,
  Kluwer, Dordrecht, 1992.

\end{thebibliography}

\end{document}